\documentclass[review]{elsarticle}

\usepackage{hyperref}
\usepackage[hyperref]{xcolor}

\hypersetup{
  colorlinks=true,
}

\usepackage{xcolor}
\usepackage{amsmath}
\usepackage{amsfonts}
\usepackage{graphicx}
\usepackage{subfigure}
\newtheorem{thm}{Theorem}[section]
\newtheorem{lem}{Lemma}[section]
\newdefinition{rmk}{Remark}[section]
\newdefinition{remark}{Remark}[section]
\newtheorem{example}{Example}[section]
\newproof{pf}{Proof}
\renewenvironment{pf}{\noindent{\textit{Proof.}  }}{\hfill $\Box$}

\DeclareMathOperator{\sgn}{sgn}
\newtheorem{prop}{Proposition}[section]

\def\RR{\mathbb{R} }
\def\NN{\mathbb{N}}
\def\EE{\mathbb{E}}
\def\cL{{\cal L}}
\def\de{{\delta}}
\def\ga{{\gamma}}
\def\Ga{{\Gamma}}
\def\be{{\beta}}
\def\th{{\theta}}
\def\var{{\rm Var}}
\def\cF{{\cal F}}
\def\bfone{{\bf 1}}
\def\vare{{\varepsilon}}

\newcommand{\al}{\alpha}

\newproof{pot1}{Proof of Theorem \ref{th:1}}
\newproof{pfL1}{Proof of Lemma \ref{e:Lop}}
\newproof{pfL2}{Proof of Lemma \ref{lem:1}}
\newproof{pfL3}{Proof of Lemma \ref{lem:2}}
\numberwithin{equation}{section}

\begin{document}

\begin{frontmatter}

\title{Local time of infinite time horizon Brownian  bridge}


\author[mymainaddress]{Yaozhong Hu}
\ead{yaozhong@ualberta.ca}

\author[mysecondaryaddress]{Yuejuan Xi\corref{mycorrespondingauthor}}
\cortext[mycorrespondingauthor]{Corresponding author}
\ead{yjx@mail.nankai.edu.cn}

\address[mymainaddress]{Department of Mathematical and Statistical Sciences, University of Alberta at Edmonton Edmonton, Alberta Canada, T6G 2G1}

\address[mysecondaryaddress]{School of Mathematical Sciences, Nankai University, Tianjin, PR China, 300071}

\begin{abstract}
We introduce an infinite time horizon Brownian bridge  which is determined  by a stochastic Langevin equation with time dependent drift coefficient. We show that this process goes to zero almost surely when the time goes to infinity and study the existence and asymptotic behavior of  its local time as well as its  H\"older continuity in time variable  and in location variable.  The main difficulty is the lack of stationarity of the process so that the powerful tools for stationary (Gaussian) processes are not applicable. We employ the Garsia-Rodemich-Rumsey inequality to get around this type of difficulty. 
\end{abstract}


\end{frontmatter}


\section{Introduction}
 The Markov bridge is widely used in statistics, probability and finance.
For example, in statistics Brownian bridge  plays an important role in the Kolmogorov-Smirnov test. In the probability theory, it is well-known that Brownian, Gamma and Bessel bridges have been developed extensively in the literature  (see e.g., \cite{bertoin1994path, emery2004parallel, pitman1982decomposition}) and references therein. By means of $h$-function,   some general Markov bridges with the SDE representation on $[0, T]$ were 
 constructed in  \cite{ccetin2016markov}. 
In finance   the $\alpha$-Brownian bridge was used in \cite{brennan1990arbitrage}  to model the arbitrage profit associated with a given futures contract.  The phenomenon of stock pinning on option expiration dates   was described  via the bridge in \cite{avellaneda2003market}. Markov  bridges  are also  employed to solve  the famous 
insider trading Kyle-Back models (cf. \cite{ back1992insider, kyle1985continuous}
and many references which follow). There are many  
other studies and  applications of Markov bridges, we  refer to e.g.
\cite{ccetin2018dynamic} and the references therein.

On a complete filtration probability space $\Lambda:=(\Omega,\mathbb P, \mathcal F_t, \mathcal F)$, a  $\al$-Brownian bridge on the interval $[0, T]$ is defined by
\begin{equation}\label{e:BB}
    \left\{
     \begin{aligned}
        &dX_t=\frac{-\al X_t}{T-t}dt+ dW_t\;, &t\in[0,T),\\
        &X_0=0, 
     \end{aligned}
     \right.
\end{equation}
where $\al>0$, $T\in(0,\infty)$, and $\{W_t\}_{t\ge0}$ is a standard  Brownian motion
on $(\Omega,\mathbb P, \mathcal F_t, \mathcal F)$. If $\al=1$, $\{X_t\}_{t\ge0}$ is the  standard Brownian bridge.  There is also  a vast literature on the property of this type of  Brownian bridge and its application, including local time and stopping time (cf. \cite{ bingham1987empirical,  ekstrom2020optimal, ekstrom2009optimal, pitman1999distribution, shepp1969explicit}).

In application, sometime we don't specify the terminal time $T$. For example,
in  the optimal portfolio and consumption problem in mathematical finance, 
sometimes it is more desirable to consider the optimal portfolio problem for an individual's life time, which
is  not a priori determined. Sometime  the terminal time $T$ we are concerned 
is    large. In this case it is convenient  to use the Brownian bridge 
over  an infinite time horizon.  
Although \cite{ccetin2018dynamic} provided the weak condition of Markov bridges with $T=\infty$, it is still in the framework of SDE containing parameter $T$. 
To the best of our knowledge, there is no reference concerning 
Brownian bridge in  the case $T=\infty$.  

To construct such an infinite time horizon Brownian bridge, we cannot simply let the drift term 
$b_T(t,x)=-\frac{\al x}{T-t}$   in  
\eqref{e:BB} go to infinite since it simply goes to zero, which yields  the Brownian motion. 
On the other hand, if the drift term is $-\al x$ then the solution  $X_t$ 
 is the famous time-homogeneous OU process,
which has a (non zero) limit as $t\rightarrow \infty$. This motivates us to require the drift 
to have the form $b(t,x)=-\al(t) x$ such that $\al(t)\rightarrow \infty$ when $t\rightarrow \infty$. 
Thus,   in this paper we deal with the general Brownian  bridge on an infinite time horizon satisfying  
the following stochastic differential equation (SDE): 
\begin{equation}\label{e:FBB}
    \left\{
     \begin{aligned}
        &dX_t=-\al(t)X_tdt+ dW_t\;, &t\in[0,\infty),\\
        &X_0=0,
     \end{aligned}
     \right.
\end{equation}
where $\al(t)$ is a deterministic  function. This equation can  be easily modified to construct more general Brownian
bridge $X_t^{a, b}$  starting at $X_0^{a, b}=a$ and 
terminating at $X_\infty^{a, b}=b$:
\begin{equation}\label{e:FBBab}
    \left\{
     \begin{aligned}
        &dX_t^{a, b}=-\al(t)(X_t^{a, b}-b)dt+ dW_t\;, &t\in[0,\infty),\\
        &X_0^{a, b}=a,
     \end{aligned}
     \right.
\end{equation}
since the solutions to the above equations \eqref{e:FBB} and \eqref{e:FBBab} are related by
\[
X_t^{a,b}=b+(a-b)e^{-\int_0^t \al(r) dr}+X_t\,.
\]

After proposing the candidate \eqref{e:FBB} for the infinite time horizon Brownian bridge  
the first task is to show we do indeed have $\displaystyle \lim_{t\rightarrow \infty} X_t=0$
almost surely. 
In the classical case of finite time horizon,  the terminal pinning point 
almost surely limit $\displaystyle\lim_{t\to T}X_t=0$   can be proved by the strong law of large numbers for Brownian motion (cf. \cite[Page 359]{karatzas1991brownian}). However, it cannot be conducted in our case because of $T=\infty$ and because of the explosive behavior of the drift term $\al(t)$.  In general, since an important feature of the process $X_t$ is that it is no longer stationary, many powerful tools effective for stationary process is no longer applicable and we need completely different probabilistic techniques.

Our strategy to prove  the bridge property of $\displaystyle\lim_{t\to \infty}X_t=0$ almost surely is first to prove   $\displaystyle\lim_{n\to \infty}X_n=0$ almost surely along positive integers $\NN$ by using some moment estimates and the Borel-Cantelli lemma. The main difficulty lays in the case of continuous time. We find that  the Garsia-Rodemich-Rumsey inequality can play a crucial role. This is done 
in Section \ref{sec:2}.  

To study a stochastic process, an important concept is its local time. For standard Brownian motions 
there are many works and here we only refer to \cite{marcusrosen}. Since our Brownian bridge is defined 
on the  half line  $\RR_+$, we are first concerned with the asymptotic behavior of 
 the local time $\cL_t^x=\int_0^t \de( X_s-x) ds$ as $t\rightarrow \infty$.  We show the existence of
the local time and we show $\displaystyle \lim_{t\rightarrow \infty}
\cL_t^x=\infty$  under some minor restrictions  on $\al(t)$  in Section \ref{sec:3}.

The H\"older continuity of the local time has been well-studied for Brownian motion and other processes (cf. \cite{MR2446295, marcusrosen} and references therein). We shall also study this property for our local time of the Brownian bridge. Again since the local time can be defined for all $t\in \RR_+$ and $x\in \RR$, it is interesting to know how the H\"older coefficient depends 
on the size of the domain. For the H\"older continuity on time variable, we shall give a positive answer.
More precisely, we shall prove that there is a finite random constant $C$, independent of $s, t, T$,  such that
$|\cL_t^x-\cL_s^x|\le   C\left[|t-s|^{1/2}    \sqrt{(T+1)\al^*(T+1) }  
+  |t-s|^{1/2} \sqrt{\log \frac{1}{|t-s|}} \ \right] $, \ $\forall\
0\le s , t\le T<\infty$. This will need some nondeterminism results for the Brownian bridge process. Compared with the result of \cite{kono} this result is sharp on any bounded interval. 
But for the H\"older continuity with respect to location parameter we can only show that its H\"older exponent can arbitrarily close to $1/2$. 
This is done in Section \ref{sec:4}.    In section \ref{sec:5}, we give some numerical experiments on a typical  sample path of  an infinite time horizon  Brownian bridge. 

The  infinite time horizon Brownian bridge process is  a special case of time-dependent modulated drift Ornstein-Uhlenbeck (OU)  processes, which is popular for modeling the evolution of interest rate (cf. \cite{ brigo2007interest, hull1990pricing}). From a theoretical point of view, it is also interesting to know that the time-dependent OU process  becomes an infinite time horizon bridge process when the drift term has   some growth rate.  

\section{Infinite time horizon Brownian  bridges}
\label{sec:2}
In this section, we  provide a sufficient condition on the drift coefficient $\al(t)$   
in Equation \eqref{e:FBB} 
so that $X_\infty=0$ (more precisely $\displaystyle \lim_{t\rightarrow\infty}X_t=0$ almost surely).  
Namely,  we provide a sufficient condition so that we have an infinite time horizon
Brownian bridge: a mean zero Gaussian process $(X_t, 0\le t\le \infty)$ such that $X_0=X_\infty=0$.  
It worths   to point out that such process is not unique as we shall see  and  from the construction 
(see Equation \eqref{eq:Sol} below), the Brownian bridge is adapted to the Brownian motion $W$. 
The main techniques  that we need are  Garsia-Rodemich-Rumsey inequality and Borel-Cantelli lemma.
First we have the following representation of the solution to \eqref{e:FBB}. 
\begin{prop}
The SDE \eqref{e:FBB} has a unique strong solution given by
\begin{equation}\label{eq:Sol}
X_t= \int_0^t\exp\left(-\int_ s^t \alpha(u)du\right)dW_s,\quad  0\leq t<\infty.
\end{equation}
\end{prop}
\begin{pf}
Since for $0\leq t<\infty$, $-\alpha(t)$ satisfies the usual Lipschitz and linear growth conditions,  the SDE \eqref{e:FBB} has a unique strong solution (see e.g. \cite[Section 5]{karatzas1991brownian}). 
By It\^o's formula, it is easy to verify that  the process $X_t$ defined by equation \eqref{eq:Sol} is the solution to \eqref{e:FBB}.
\end{pf}

Next, we    define 
\begin{equation}\label{e:X}
    X_t:=\left\{
     \begin{aligned}
        & \int_0^t\exp\left(-\int_s^t\alpha(u)du\right)dW_s,&t<\infty,\\
        &0,&t=\infty\,. 
     \end{aligned}
     \right.
\end{equation}
For convenience of notations,  throughout this paper we shall use $C$ to denote a generic finite 
positive constant, whose value may be different in different appearances. 
When we need to stress the dependence of a constant on $\ga$, we use  $C_\ga$, whose values  may also be different in different appearances.   
 For a function $\alpha(t):[0,\infty)\to\mathbb R $ we introduce the following notation:
 \begin{equation}
\al^*(t):=\sup_{0\le r\le t} \al(r)\,, \quad 0\le t<\infty\,. 
 \end{equation}
\begin{thm}\label{th:1}
  Let  the time-dependent coefficient  function $\alpha(t):[0,\infty)\to\mathbb R^+$ be   continuously differentiable  and satisfy the following growth condition:
\begin{enumerate}[(i)]
  \item There are  $\ga\in [0, 1/2]$ and    $\beta>0$  such that $\frac{(\al^*(t+1))^{ 2\ga}}{\al(t) }\le C  t^{-\beta} $ for all $0\le t<\infty$.
  \item There is a constant $C$ such that $\sup_{t\ge 0}|\frac{\al^\prime(t)}{\al(t)}|\le C$.  
\end{enumerate}
 Then the process $\{X_t\}_{t\geq0}$ defined by \eqref{e:X} is a centered Gaussian process with almost surely continuous sample paths on the closed infinite time horizon $[0, \infty]$.
\end{thm}

\begin{remark} 
\begin{enumerate}
\item[(1)]Condition (i) implies that 
\begin{equation}
\al(t)\geq C_\be t^{\beta}\,, \quad \forall \ 0\leq t<\infty\,. 
\end{equation}
\item[(2)] The condition (i) in Theorem \ref{th:1} can be replaced by
``there are  $t_0>0$,  $\ga\in [0, 1/2]$,   $\beta>0$,  and $K>0$ 
  such that $\frac{(\al^*(t+K))^{ 2\ga}}{\al(t) }\le C  t^{-\beta} $ for all $t_0\le t<\infty$". 
\item[(3)] The condition that $\al(t)\rightarrow \infty$ seems necessary. For example if $\al(t)=\al$ is 
a constant function, the process $X_t$ will be a usual OU process, which converges 
to a nondegenerate Gaussian random variable. 
\end{enumerate}
\end{remark} 

To prove this theorem, we need  the following three lemmas, whose proofs are given in  \ref{A:A}
\begin{lem}\label{e:Lop}
Assume that the function $g:[0,\infty)\to \RR^+$ is  continuous and satisfies $\lim_{x\rightarrow +\infty}  g(x) =+\infty$. If there is a positive constant $C$
such that $ \sup_{x>0}  |\frac{g'(x)}{g(x)}|\le C$,   then for any $\kappa>0$, there are positive constants $x_0$, $\bar C_{\kappa,g}$, and $C_{\kappa, g}$,
independent of $x\in \RR_+$,  such that for all $x>x_0$,
\begin{equation}
\bar C_{\kappa,g}\le\frac{g(x) \int_0^x\exp\left(\kappa \int_0^s g(u)du\right)ds}{\exp\left( \kappa \int_0^x g(u)du\right)  }\leq   C_{\kappa, g}\,,     \label{e.2.9} 
\end{equation}
\end{lem}

\begin{lem}\label{lem:1}
 Let $X$ be the solution to \eqref{e:FBB} and  $\al(t)$ satisfy the conditions in Theorem \ref{th:1}. Then for any $\ga\in [0,1/2]$ and  for any  $t_1, t_2>0$,
\begin{align}
\sigma^2_{t_1, t_2}&:=\mathbb E|X_{t_2}-X_{t_1}|^2\nonumber\\
&\le C_{ \ga}|{t_2}-{t_1}|^{2\ga}\left[ \left(\al^*(t_1\vee t_2)\right)^{2\ga}   (\al(t_1\wedge t_2))^{-1} +
(\al(t_1\vee t_2))^{ 2\ga-1} \right] \,. \label{e.2.6} 
\end{align}
\end{lem}

\begin{example} If there is a $t_0>0$ such that for all $t>t_0$, $\al(t)\ge Ct^{\be}$ for some $\be>0$ and $\al(t)\le Ct^p$ for some $\be \le p<\infty$, and if $t_0<t_1<t_2<t_1+1$,  then
\begin{eqnarray*}
\sigma^2_{t_1,t_2} 
&\leq&  C_{  \ga} 
 \left|t_2-t_1\right|^{2\ga} \left[\left(  t_1+1 \right)^{2\ga p}    t_1^{-\be}  +    (  t_1 )^{-\beta(1-2\ga)}\right]   \\
 &\leq&  C_{ \ga} 
 \left|t_2-t_1\right|^{2\ga}  t_1^{-\rho}   
\end{eqnarray*}
for some $\rho\in (0, \be)$ when $\ga>0$ is sufficiently small. We can also assume that $\al(t)>e^{pt}$ and $\al(t)\le e^{qt}$ for some positive constants $p$ and $q$. 
\end{example}

\begin{lem}\label{lem:2} Let $\beta>0$, $\ga\in [0, 1/2]$.
 Fix an arbitrary  positive number $N>1$. Choose a   positive integer $m$ such that
 $m\ga>2$, $2N/m<\beta$, $2\ga<1-2N/(\beta m)$.  Assume that the conditions   in Theorem \ref{th:1} are satisfied. Then, there is a random constant 
 $R_{N,m}$  (independent of $k$) such that 
for any integer $k\geq1$, $t_1, t_2 \in[k,k+1]$, the following inequality holds.
\begin{equation}
|X_{t_2}-X_{t_1}|\le R_{N,m} C_\ga k ^{-\beta\left[1-2\ga-2N/(\beta m)\right]/2} \,. \label{e.2.10} 
\end{equation}
\end{lem}

Now we are ready to prove Theorem \ref{th:1}.

\begin{pot1}
  First, we show 
$\lim_{n\rightarrow \infty}  X_n=0$  along the positive integers.   For any real number $\epsilon>0$, any positive integer $q$ with $\be q>1$, by the Chebyshev's inequality and then by Lemma  \ref{e:Lop}, we have
\begin{align*}
\sum_{k=1}^\infty \mathbb P(|X_k|>\epsilon)&\leq\sum_{k=1}^\infty\frac{\mathbb E|X_k|^{2q}}{\epsilon^{2q}}\\
&\le \sum_{k=1}^\infty\frac{M_q\left(\int_0^k\exp\left(-\int_s^k2\alpha(u)du\right)ds\right)^{q}}{\epsilon^{2q}}\\
&\leq\sum_{k=1}^\infty\frac{C_q }{\epsilon^{2q}\left(\al(k)\right)^q} \leq\sum_{k=1}^\infty\frac{C_q }{\epsilon^{2q}k^{\be q}} <\infty,
\end{align*}
where $M_p=\frac{(2q)!}{2^{q}(q)!}$ . Then the Borel-Cantelli lemma implies $\mathbb P(|X_k|>\epsilon ~\textup{i.o.})=0$. Since $\epsilon$ is arbitrary, then $X_k$ converges to zero almost surely as $k\to\infty$. 

Next, we will show the almost sure convergence    result for   continuous time, namely, we need to show $X_t$ converges to zero almost surely as $t\to\infty$. Clearly, we have  
\begin{equation}\label{eq:C}
|X_t|\leq|X_t-X_k|+|X_k|,
\end{equation}
where $k=\lfloor t\rfloor$ is the biggest integer less than or equal to the  real number $t$. The second term in \eqref{eq:C} converges to zero almost surely as $k=\lfloor t\rfloor\to \infty$ from the previous argument.
Lemma \ref{lem:2} can be used to show  that the first term in equation
\eqref{eq:C} also converges to zero almost surely as $t\to\infty$.\qed
\end{pot1}

\section{Local time of infinite time horizon Brownian bridges}
\label{sec:3}
 In this section, we  study the local time for the infinite 
 time horizon Brownian  bridges. 
The local time of the one-dimensional infinite time horizon 
Brownian  bridge $\{X_t\}_{t\ge0}$ at level $x$ is defined by
\begin{equation}
\mathcal L_t^x=\int_0^t\delta(X_r-x)dr,
\end{equation}
which is the limit (in $L^2$ if the limit exists) of the approximating  local time process   defined by
  \begin{equation}
\mathcal L_{t,\varepsilon}^x =\lim_{\varepsilon\to0+}\int_0^tp_\varepsilon(X_r-x)dr,
\end{equation}
where $\delta(\cdot)$ is the Dirac delta function at zero, $p_\varepsilon(x)=\frac{1}{\sqrt{2\pi\varepsilon}}\exp\left(-\frac{x^2}{2\varepsilon}\right)$ is the  heat kernel  with $\varepsilon>0$.
To study the limit of $\mathcal L_{t,\varepsilon}^x$ we use the  following representation  for the   heat kernel: 
\begin{equation}
p_\varepsilon(x)=\frac{1}{2\pi}\int_\RR \exp\left( i\xi x-\frac{\varepsilon\xi^2}{2}\right)d\xi,
\label{e.3.3}
\end{equation}
where $i^2=-1$.
\begin{lem}\label{lem:Cauchy}
For any $x\in \RR$, $0< t<\infty$,  we have 
\begin{equation*}
\lim_{\varepsilon,  \theta\to0+}\mathbb E|\mathcal L_{t,\varepsilon}^x-\mathcal L_{t,\theta}^x|^2=0.
\end{equation*}
\end{lem}
\begin{pf} 
To simplify notation, and without loss of generality we will assume  $x=0$ in the following argument for the existence of the limit. We set $\mathcal L_t=\mathcal L_t^0$ and $\mathcal L_{t,\varepsilon}=\mathcal L_{t,\varepsilon}^0$.
Since $\mathbb E|\mathcal L_{t,\varepsilon}-\mathcal L_{t,\theta}|^2=\mathbb E(\mathcal L_{t,\varepsilon})^2+\mathbb E(\mathcal L_{t,\theta})^2-2\mathbb E(\mathcal L_{t,\varepsilon}\mathcal L_{t,\theta})$, 
we only need to compute $\mathbb E(\mathcal L_{t,\varepsilon}\mathcal L_{t,\theta})$.  We use the expression 
\eqref{e.3.3} for this computation: 
\begin{align*}
\mathbb E(\mathcal L_{t,\varepsilon}\mathcal L_{t,\theta})&=  \frac{1}{4\pi^2}\int_0^t\int_0^t  \int_{ \RR^2}
\mathbb E \left\{\exp\left( i(\xi X_r+\eta X_s) -\frac{\varepsilon\xi^2+\theta\eta^2) }{2}\right)\right\}d\xi d\eta dsdr\\
&= \frac{1}{4\pi^2}\int_0^t\int_0^t  \int_{ \RR^2}
  \left\{\exp\left(-\frac12   E \left[\xi X_r+\eta X_s\right]^2 -\frac{\varepsilon\xi^2+\theta\eta^2) }{2}\right)\right\}d\xi d\eta  dsdr\\
  &= \frac{1}{2\pi}\int_0^t\int_0^t  \left[ {\rm det} (A_{\varepsilon,\theta}(s,r)) \right]^{-1/2}  dsdr\\
  &= \frac{1}{ \pi}\int_0^t\int_0^s  \left[ {\rm det} (A_{\varepsilon,\theta}(s,r)) \right]^{-1/2}  drds,
\end{align*}
where
\begin{equation}
A_{\varepsilon,\theta}(s,r)= \begin{pmatrix}
\mathbb  E(X_s^2) + \theta&\mathbb E(X_rX_s)\\
\mathbb E(X_rX_s)& \mathbb E(X_r^2)+\varepsilon\\
\end{pmatrix}.
\end{equation}  
Now we need to compute the above determinant. If $0< r<s< t<\infty$,  then 
\begin{align*}
\rm det (A_{0,0}(s,r))&=\EE X_s^2\EE X_r^2-(\EE X_r X_s)^2\\
&=\exp\left\{ -2\int_0^r\al(u)du-2\int_0^s\al(u)du\right\}  \int_0^r\exp\left(2\int_0^u\al(v)dv \right)du \\
&  \int_r^s\exp\left(2\int_0^u\al(v)dv \right)du \\
&\ge  \exp\left\{  -2\int_0^s\al(u)du\right\}  \int_0^r\exp\left(2\int_0^u\al(v)dv \right)du\cdot  (s-r) \\ 
&\ge  C_t r  (s-r)  , 
\end{align*}
for some constant $C_t>0$,  where the last inequality follows from the fact  that a continuous function  attains  the minimum and maximum on the bounded interval $[0,t]$.  It is easy to see that 
for any $\th, \vare>0$,  $\rm det (A_{ \vare, \th}(s,r))\ge \rm det (A_{0,0}(s,r))$
and 
\[
\lim_{\vare, \th\rightarrow 0+} \rm det (A_{ \vare, \th}(s,r))=\rm det (A_{0,0}(s,r))
\quad \hbox{for any $0< r<s<t$}\,.
\]
On the other hand for any $0\le p<2$, 
\begin{align*}
\int_0^t\int_r^t[\rm det (A_{0,0}(s,r))]^{-p/2}dsdr&\le C_t \int_0^t\int_r^t r^{-p/2}(s-r)^{-p/2}dsdr\\
&\le C_t \kappa^2t^{2-p} <\infty,
\end{align*}
for some constant $\kappa$, where the last equality above follows from \citet[Lemma 9.1]{yaozhong2017analysis}.
By the Lebesgue's dominated convergence theorem, we have 
\begin{equation*}
\lim_{\varepsilon\to0+,\theta\to0+}\mathbb E(\mathcal L_{t,\varepsilon}\mathcal L_{t,\theta})=
\frac{1}{ \pi}\int_0^t\int_0^s  \left[ {\rm det} (A_{0,0}(s,r)) \right]^{-1/2}  drds<\infty\,. 
\end{equation*}
 Since $\varepsilon$ and $\theta$ are arbitrary, we can obtain that 
\begin{equation*}
\lim_{\varepsilon\to0+}\mathbb E (\mathcal L_{t,\varepsilon})^2=\mathbb E (L_t)^2 \textup{    and }  \lim_{\theta\to0+}\mathbb E (\mathcal L_{t,\theta})^2=\mathbb E (L_t)^2.
\end{equation*}
Then the desired result is proved.
\end{pf} 

Lemma \ref{lem:Cauchy} immediately implies the following theorem.
\begin{thm}
For any $0<t<\infty$ and any $x\in \RR$,  we have
\begin{equation}
\lim_{\varepsilon\to0+}\mathbb E|\mathcal L_{t,\varepsilon}^x-\mathcal L_t^x|^2=0.
\end{equation}
\end{thm}
Unlike the classical Brownian bridge, ours is defined for all $t\ge 0$. Thus, the local time is also 
well-defined for all time $t\ge 0$. It is then arisen an interesting question: does the local time 
$\cL_t$ have a limit as $t\rightarrow \infty$?  Intuitively, $\de(X_s)$ is nonnegative so that $\int_0^t \de(X_s)ds$ is an increasing (in time variable $t$)  stochastic process. So the limit of $\int_0^t \de(X_s)ds$  as $t\rightarrow \infty$ should exist as a finite or infinite variable.
In the  next theorem, we show that  the local time process $\mathcal L_t$ goes to infinite when $t$ tends to infinity.

\begin{thm}Let $X_t$ be the Brownian bridge satisfying the conditions of Theorem \ref{th:1} and 
assume that there is a $\rho>1$ such that $\al$ satisfies  
\begin{equation}
\alpha(t)\ge Ct^\rho\,, \quad \forall \ t\ge t_0\quad \hbox{for some $  t_0>0$}\,. \label{e.3.6}
\end{equation} 
Then
\begin{equation}
\lim_{t\to\infty}\mathcal L_t=\infty\quad  \hbox{almost surely}\,.
\end{equation}
\end{thm}

\begin{pf}  
By the It\^o-Tanaka formula (see e.g. \citet{revuzyor}), we have
\[
|X_t-x|=|X_0-x|+\int_0^t \sgn (X_r-x)dB_r-\int_0^t \sgn(X_r-x)
\alpha(r)  X_rdr+\int_0^t \delta(X_r-x) dr\,. 
\]
Or 
\begin{equation*}
\int_0^t \delta(X_r-x) dr=|X_t-x|-|X_0-x|-\int_0^t \sgn (X_r-x)dB_r+\int_0^t \sgn(X_r-x)
\alpha(r)  X_rdr\,, 
\end{equation*} 
where $\sgn (x)$ denotes the sign of the real number $x$. 
Letting $x=0$ yields 
\begin{equation}
\cL_t=\int_0^t \delta(X_r) dr=|X_t| -\int_0^t \sgn (X_r)dB_r+\int_0^t   
\alpha(r)  |X_r|dr\,. \label{e.3.8} 
\end{equation} 
It is clear that 
$(\int_0^t \sgn (X_r)dB_r, T\ge 0)$ is a Brownian motion. So, 
for any $\nu >1/2$, there is a (random) constant $C_\nu$ such that 
\[
\left|\int_0^t \sgn (X_r)dB_r\right|\le C_\nu t^\nu \,,   \quad \forall \ t\ge 0\,.
\]
Next, we want to show that  there is  a $\mu>1/2$ such that
\[
\liminf_{t\rightarrow \infty} \frac{1}{t^\mu}\int_1^t   
\alpha(r)  |X_r|dr>0\,. 
\]
For any $M>0$ sufficiently large,   any  $p\in(0,1)$, any $\mu>0$,   and any $t>1$,   by the Chebyshev's 
inequality  we have 
\begin{eqnarray*}
P(\frac{1}{t^\mu}\int_1^t   
\alpha(r)  |X_r|dr\le M)
&=& P( M^p t^{p\mu} \left(\int_1^t   
\alpha(r)  |X_r|dr \right)^{-p}  \ge 1)\\
&\le&   M^p t^{p\mu} \EE \left(\int_1 ^t   
\alpha(r)  |X_r|dr \right)^{-p} \\
&=&  M^p t^{p\mu} (t-1)^{-p}  \EE \left(\frac{1}{t-1} \int_1^t   
\alpha(r)  |X_r|dr \right)^{-p}  
\end{eqnarray*}
Since when $p\in(0,1)$,
$\phi(x)=x^{-p}$ is a convex function  an application of the Jensen's inequality 
to $\frac{1}{t-1} \int_1^t f(r)   dr$ yields
\begin{eqnarray}
P(\frac{1}{t^\mu}\int_1^t   
\alpha(r)  |X_r|dr\le M)
&\le&  M^p t^{p\mu} (t-1)^{-p}   \EE \left(\frac{1}{t-1} \int_1^t   
\alpha(r) ^{-p}  |X_r|^{-p} dr \right)  \nonumber \\ 
&=& M^p t^{p\mu} (t-1)^{-p-1}    \int_1^t   
\alpha(r) ^{-p}  \EE |X_r|^{-p} dr  \,.  \label{e.3.9} 
\end{eqnarray}
On the other hand, by Lemma \ref{e:Lop}  we have 
 \begin{eqnarray*}
\sigma_r^2:=\EE(X_r^2)
&=& \exp\left(-2\int_0^r\alpha(s)ds\right)\int_0^r\exp\left(2\int_0^s\alpha(u)du \right)ds \ge  C/\alpha(r) .
\end{eqnarray*}
This implies that for any $p\in (0,1)$, 
\begin{eqnarray*}
\EE(|X_r|^{-p} )&=&\int_\RR |x|^{-p}\frac{1}{\sqrt{2\pi\sigma_r^2}}\exp\left( -\frac{x^2}{2\sigma_r^2}\right)dx\\
                        &=&\sigma_r^{-p}\int_\RR |y|^{-p}\frac{1}{\sqrt{2\pi}}\exp\left( -\frac{y^2}{2}\right)dy\\
                       &\le&   C_p \left(\alpha(r)\right) ^{p/2}  \,. 
\end{eqnarray*}
Substituting the estimate into \eqref{e.3.9} and using the 
condition \eqref{e.3.6}, we obtain for $t$ sufficiently large 
(e.g. $t\ge 2$)
\begin{eqnarray}
P(\frac{1}{t^\mu}\int_1^t   
\alpha(r)  |X_r|dr\le M)  
&\le& M^p t^{p\mu} (t-1)^{-p-1}    \int_1^t   
\alpha(r) ^{-p/2}    dr  \nonumber \\  
&\le  &  C_p M^p t^{p\mu}t^{-p-1}     \int_1^t   
 r^{-\rho p/2}    dr  \nonumber  \\ 
 &\le  &  C_p M^p t^{p\mu} t^{-p-1}    \left[t^{-\frac{\rho p}{2}+1}+1\right]  \nonumber  \\ 
  &\le  &
 \begin{cases}
  C_p   t^{p\mu -p-\frac{\rho p}{2}}         &\qquad \hbox{if $\rho <2$}\\  
  C_p  M^p  t^{p\mu  -p-1 }       &\qquad \hbox{if $\rho >2$} 
 \end{cases}\label{e.3.10} 
\end{eqnarray}  
for $p$ sufficiently close to $1$ since  when $\rho>2$,  $t^{-\frac{\rho p}{2}+1}$ is bounded for $t\ge 2$. 
When $p=1$, the two exponents in \eqref{e.3.10} are
$p\mu -p-\frac{\rho p}{2}=\mu-\frac{\rho}{2}-1$ and $p\mu-p-1=\mu-2$.  
From these computations, we see clearly that when  
\begin{equation}
\mu <\min(\rho/2, 1)\,, 
\end{equation}
we can choose   $p$ sufficiently close  to $1$ so that both exponents in 
\eqref{e.3.10} will be less than $-1$. Namely,  we can find   an  $\ell >1$ (with an appropriate choice of $p$ close to $1$
in \eqref{e.3.10}) such that
\begin{equation*}
P\left(\frac{1}{t^\mu}\int_1^t   
\alpha(r)  |X_r|dr\le M\right)\le C t^{-\ell }\,. 
\end{equation*}
This implies  
\begin{equation*}
\sum_{n=1}^\infty P\left(\frac{1}{n^\mu}\int_1^ n   
\alpha(r)  |X_r|dr\le M\right)<\infty\,. 
\end{equation*} 
By the Borel-Cantelli lemma we see that 
\begin{equation}
\lim_{n\rightarrow \infty} \frac{1}{n^\mu}\int_1^ n   
\alpha(r)  |X_r|dr=\infty\,. \label{e.3.11} 
\end{equation}
Dividing both sides of \eqref{e.3.8} by $n^\mu$ we get
\begin{equation}
\frac{1}{n^\mu}\cL_n=\frac{1}{n^\mu}\int_0^t   
 \alpha(r)  |X_r|dr+I_{1,n}+I_{2,n}
\,,  \label{e.3.12}
\end{equation} 
where $I_{1,n}=\frac{1}{n^\mu}|X_n| $ and $I_{2,n}=-\frac{1}{n^\mu}\int_0^n \sgn (X_r)dB_r$. 
Since $|X_n|\stackrel{\rm a.s.} \rightarrow 0$, we see that $I_{1,n}
\stackrel{\rm a.s.} \rightarrow 0$.  Since we assume $\rho>2$ and 
$\nu>1/2$ is arbitrary, we can choose $\mu>\nu$. Thus,  we also have $I_{2,n}\stackrel{\rm a.s.} \rightarrow 0$.
Combining the above results with \eqref{e.3.11}-\eqref{e.3.12},  we see that 
\begin{equation}
\lim_{n\rightarrow \infty} 
\frac{1}{n^\mu}\cL_n  =\infty  \quad {\rm almost\ surely} 
\end{equation} 
which in turn  implies
\begin{equation}
\lim_{n\rightarrow \infty} 
 \int_0^n \delta(X_r) dr  =\infty\quad {\rm almost\  surely} \,. 
\end{equation} 
Since $\cL_t $ is increasing on $t$ almost surely
(it is the limit of a sequence of the approximating local times processes $\cL_{t, \vare}$
which is obviously  increasing in $t\ge 0$), we see that
\begin{equation}
\lim_{t\rightarrow \infty} 
 \cL_t  =\infty\,. 
\end{equation} 
This completes the proof of the theorem. 
\end{pf} 
\begin{remark} The above argument show that $\cL_n\ge C n^{\nu}\,, \forall \ n\in \NN$ for any $\nu<\rho/2$. 
It is natural to conjecture that $\cL_t\ge C t^{\rho/2}$  for all positive $t$ sufficiently  large.
\end{remark} 

\section{H\"older continuity of local time}\label{sec:4} 
The H\"older continuity of the local time 
 of a stochastic process is always an important topic in the probability theory. 
 Since the local time $\cL_t^x$ of the infinite time horizon Brownian bridges 
 depends on two parameters: time parameter $t$ and location parameter $x$ 
we shall study
 in this section the H\"older continuity of  $\cL_t^x$ with respect to $t$ and with respect to $x$ separately.  
We know that the H\"older constant usually depends on the (bounded) domain we are working on. 
Since our Brownian bridges are defined on the whole half line, we are interested in 
the problem how the H\"older constant depends on the size $T$  of the domain $[0, T]$. 
We have a positive answer for this  problem with respect to the time parameter. But it 
seems  hard to work on the 
location parameter $x$.  Our result  on the H\"older continuity of the local time with respect to
the time is more precise than that  for location parameter $x$.

We shall use the following results which is analogous to the nondeteministic results on our Brownian bridge process $X_t$. But we also need a  upper 
bound estimate.   
We  assume that the function $\al(\cdot)$ is measurable and positive. 
\begin{lem}\label{lem:5}
Let $X_t$ be the solution to the equation \eqref{e:FBB} and let $p$ be a positive integer. 
For  $u=(u_1, \cdots, u_p)$ with $0\le u_1< u_2< \cdots  <u_p<\infty$, denote 
\begin{equation}
A_p(u)=\left(a_{ij}(u)\right)_{1\le i,j\le p}\,,\quad \hbox{with}\quad 
a_{ij}(u):=\EE(X_{u_i}X_{u_j})\,. \label{e.4.1} 
\end{equation}
Then   
\begin{equation} 
u_1(u_2-u_1)\cdots (u_p-u_{p-1})\exp\left\{-2\al^*(u_p) u_p\right\}\le 
 \det(A)\le u_1(u_2-u_1)\cdots (u_p-u_{p-1})\,. \label{e.4.2} 
\end{equation}
\end{lem}
We give  a proof of this  lemma    in \ref{B}. Now we  state and prove our first main result of this section on the H\"older continuity of the local time with respect to the time variable. 
\begin{thm}Fix an arbitrary  $x\in \RR$.   There exists  a (random) constant $C $ independent of $s, t,   T\in \RR_+$ 
such that 
\begin{equation}
|\cL_t^x-\cL_s^x|\le   C\left[|t-s|^{1/2}    \sqrt{(T+1)\al^*(T+1) }  
+  |t-s|^{1/2} \sqrt{\log \frac{1}{|t-s|}} \ \right]     
\end{equation}
for all $   \ 0\le  s, t\le T<\infty\,, \quad |t-s|<1$.
\end{thm} 
\begin{pf}
For  any  positive integer $p$ we first  compute the following moment: 
\begin{eqnarray*}
\EE |\cL_t^x-\cL_s^x|^p
&=& \EE |\int_s^t \de(X_u-x) du|^p\\
&=&  \int_{[s, t] ^p} \EE \left[\de(X_{u_1}-x)\cdots \de(X_{u_p}-x)\right]du\,,
\end{eqnarray*}
where $du=du_1\cdots du_p$ and where we use the Dirac function notation directly. We shall also use the formal expression $\de(x)=\frac{1}{2\pi} \int_\RR e^{i x\xi}d\xi$,  which can be justified
easily by a limiting argument through  \eqref{e.3.3}.  It is well-known that 
$\EE e^X=e^{\frac12 \EE(X^2)}$ for any mean zero Gaussian $X$.  Thus,  we have
\begin{eqnarray*}
\EE |\cL_t^x-\cL_s^x|^p 
&=&  \frac{1}{(2\pi)^p} \int_{[s, t] ^p} \int_{\RR^p}
\EE \exp\left[ -ix\sum_{k=1}^p\xi_k+i \sum_{k=1}^pX_{u_k}\xi_k\right]d\xi du\\
&=&  \frac{1}{(2\pi)^p} \int_{[s, t] ^p} \int_{\RR^p}
 \exp\left[ -ix\sum_{k=1}^p\xi_k-\frac12\EE(\sum_{k=1}^pX_{u_k}\xi_k)^2\right]d\xi du\\
&=& \frac{1}{(2\pi)^p} \int_{[s, t] ^p} \int_{\RR^p}
 \exp\left[ -ix\sum_{k=1}^p\xi_k-\frac12 \xi^\top A_p(u)\xi\right]d\xi du\,,
\end{eqnarray*}
where $d\xi=d\xi_1\cdots d\xi_p$ and $A_p(u)$ is defined by \eqref{e.4.1}. 
Integrating $d\xi$ gives 
\begin{eqnarray}
\EE |\cL_t^x-\cL_s^x|^p  
&=& \frac{1}{(\sqrt{2\pi})^p} \int_{[s, t] ^p}  \left[\det(A_p(u))\right]^{-1/2} 
 \exp\left[   -\frac {   \bfone ^\top A_p^{-1}(u)\bfone }{2} x^2 \right]  du\nonumber\\  
 &\le& \frac{1}{(\sqrt{2\pi})^p} \int_{[s, t] ^p}  \left[\det(A_p(u))\right]^{-1/2} 
   du\nonumber\\ 
 &\le& \frac{p!}{(\sqrt{2\pi})^p} \int_{s\le u_1<\cdots<u_p\le t}  \left[\det(A_p(u))\right]^{-1/2}   du\,, \label{e.3.20} 
\end{eqnarray}
where $\bfone=(1, \cdots, 1)^\top$ is the $p$-dimensional column vector whose elements are all equal to $1$
and $ A_p^{-1}(u)$ denotes the inverse matrix of $ A_p (u)$, which does exist by the first inequality in
\eqref{e.4.2}.  
 Substituting the first inequality of  \eqref{e.4.2}  to \eqref{e.3.20}, using   \cite[Lemma 4.5]{HHNT}
to integrate  $u_1, \cdots, u_{p-1}$, and denoting $I_s^t=\{s< u_1< \cdots<u_p< t\}$, then we have 
\begin{align}
\EE |\cL_t^x-\cL_s^x|^p 
&\le  \frac{p!}{(\sqrt{2\pi})^p} \int_{I_{s}^t} \left(u_1(u_2-u_1)\cdots (u_p-u_{p-1})\right)^{-\frac{1}{2}} \exp\left\{\al^*(u_p) u_p\right\}  du \nonumber\\
&\le  \frac{p!}{(\sqrt{2\pi})^p} \int_{I_{s}^t} \left((u_1-s)(u_2-u_1)\cdots (u_p-u_{p-1})\right)^{-\frac{1}{2}} \exp\left\{\al^*(u_p) u_p\right\}  du \nonumber\\
&\le \frac{p![\Ga(1/2)]^{p-1}}{(\sqrt{2\pi})^p\Ga(\frac{p-1}{2} +1)} \int_s^t
(u_p-s)^{\frac{p-2}{2}} \exp\left\{\al_*(u_p) u_p\right\}  du_p\nonumber\\
&\le \frac{C_pp!}{\Gamma(\frac{p+1}{2})} (t-s)^{\frac{p }{2}} \exp(t\al^*(t))\,.
\label{e.3.25}
\end{align}
Let us denote
\[
\rho(t)=t^{1/2}\,,\quad \Psi(t)=\exp\left(\mu t^2\right)-1\,, \quad t  \ge 0\  
\]
 for $\mu>0$.  Thus $\rho'(t)=t^{-1/2}/2$ and $
 \Psi^{-1}(t)=\frac{1}{\sqrt{\mu}} \sqrt{\log (1+t)}$. From \eqref{e.3.25} we have  for any $T>0$
\begin{eqnarray}
&&\EE \left\{ \int_0^T \int_0^T \Psi\left(\frac{\left|L_t^ x-L_s^x\right|
}{\rho(|t-s|)}\right)dsdt\right\}\le 
\sum_{p=1}^\infty \frac{\mu^p }{p!}   \int_0^T \int_0^T  
  \frac{\EE \left|L_t^ x-L_s^x\right|^{2p}
}{ |t-s| ^{p } }  dsdt \nonumber\\
&&\qquad \le \sum_{p=1}^\infty \frac{\mu^p }{p!}  \frac{C (2p)!}{\Ga(p+\frac12)}
\exp(T\al^*(T)) \int_0^T \int_0^T   dsdt\nonumber \\
&&\quad \le C_\mu T^2\exp(T\al^*(T))\label{e.3.26}
\end{eqnarray}
when $\mu$ is sufficiently small,  where we used the Stirling formula 
$\Ga(m+1)\approx  \sqrt{2\pi}m^{m+1/2}e^{-m}$ and $m!=\Ga(m+1)$.
Denote
\[
B=\sum_{n=1}^\infty \frac{1}{n^4}     \left\{ \int_0^n\int_0^n \Psi\left(\frac{\left|L_t^ x-L_s^x\right|
}{\rho(|t-s|)}\right)dsdt \exp\left(-n\al^*(n) \right)\right\}\,.
\]
Then from \eqref{e.3.26} it follows 
\[
\EE (B)\le C\sum_{n=1}^\infty \frac{1}{n^2}<\infty\,. 
\]
This means that $B$ is almost surely finite.
Since each term in the definition of $B$ is positive we know that     each  summand  
of $B$ is less  than  or equal to $B$. Thus, 
\[ 
\int_0^n\int_0^n \Psi\left(\frac{\left|L_t^ x-L_s^x\right|
}{\rho(|t-s|)}\right)dsdt\le B n^4 \exp\left\{  n\al^*(n)\right\}\,. 
\]
Since $\int_0^T\int_0^T \Psi\left(\frac{\left|L_t^ x-L_s^x\right|
}{\rho(|t-s|)}\right)dsdt$  is increase in $T$, we have
\begin{equation}
\int_0^T\int_0^T \Psi\left(\frac{\left|L_t^ x-L_s^x\right|
}{\rho(|t-s|)}\right)dsdt\le B  (T+1) ^4 \exp\left\{  (T+1)\al^*(T+1)/2\right\}\,,\quad \forall \ T\ge 0\,. 
\end{equation}
We can take $B\ge 1$.  By  the Garsia-Rodemich-Rumsey inequality (see e.g. \cite[Theorem 2.1]{yaozhong2017analysis})  we have
\begin{eqnarray}
|L_t^ x-L_s^x|
&\le & \frac{4}{\sqrt{\mu}}\int_0^{|t-s|} u^{-1/2} \sqrt{\log (1+ \frac{
4B(T+1)^4 \exp((T+1)\al^*(T+1))}{u^2} )}  du \nonumber\\
&\le & C\int_0^{|t-s|} u^{-1/2} \sqrt{\log (   
5B(T+1)^4 \exp((T+1)\al^*(T+1)) )}du \nonumber\\
&&\qquad  + C\int_0^{|t-s|} u^{-1/2}\sqrt{\log (1/u)}   du \nonumber\\
&\le& C|t-s|^{1/2}    \sqrt{(T+1)\al^*(T+1)+C\log (T+1)+C}\nonumber\\
&&\qquad + C\int_0^{|t-s|} u^{-1/2}\sqrt{\log (1/u)}   du\nonumber\\
&\le& C\left[|t-s|^{1/2}    \sqrt{(T+1)\al^*(T+1) }  
+  |t-s|^{1/2} \sqrt{\log \frac{1}{|t-s|}}\; \right]    
\end{eqnarray}
for some (random) constant $C$, independent $s,t, T $ when  $|t-s|<1$.
This shows the theorem. 
\end{pf}

Now we turn to the  H\"older continuity in $x$ of $\cL_t^x$. From the well-known results
 of  \citet[Equation (2.215)]{marcusrosen}  (see also Remark \ref{r.3.3} at the end of this section)
 it follows that 
 when $X_t$ is the   standard Brownian motion
 $(\cL_t^x \,, x\in \RR)$ is H\"older continuous of exponent $\be$ for any $\be<1/2$. 

\begin{thm}
 Fix  any positive real number $R>0$ and any $t>0$. For any $\al\in (0, 1/2)$, 
there is a positive random constant $C_{\al, t, R}$ such that
\begin{equation}
|\cL_t^x-\cL_t^y|\le C_{\al, t, R}|x-y|^{\al}\,,\quad \forall \ x, y\in [-R, R]\,. 
\end{equation}
\end{thm}  

\begin{pf}
Let $p$ be a   positive integer. We consider the  
  moment:  
\begin{eqnarray*}
\EE |\cL_t^x-\cL_t^y|^p
&=& \EE |\int_0^t \left[\de(X_u-x) -\de(X_u-y)\right] du|^p\\
&=& \int_{[0, t] ^p} \EE \left\{ \prod_{k=1}^p\left[ \de(X_{u_k}-x)-\de(X_{u_k}-y) \right]\right\}du\,,
\end{eqnarray*}
where $du=du_1\cdots du_p$. Using the formal expression for the Dirac function
 and using the notations  as in the proof of the previous theorem, we have 
\begin{eqnarray*}
\EE |\cL_t^x-\cL_t^y|^p 
&=&  \frac{1}{(2\pi)^p} \int_{[0, t] ^p} \int_{\RR^p}
\EE   \exp\left\{i \sum_{k=1}^p X_{u_k}\xi_k \right\}\prod_{k=1}^p\left[ e^{-ix\xi_k}- 
e^{-iy\xi_k}\right]   d\xi du\\
&=&  \frac{1}{(2\pi)^p} \int_{[0, t] ^p} \int_{\RR^p}
   \exp\left\{-\frac12 \xi^\top A_p(u) \xi\right\} \prod_{k=1}^p\left[ e^{-ix\xi_k}- 
e^{-iy\xi_k}\right]    d\xi du  \\
&=&  \frac{1}{(2\pi)^{p/2}} \int_{[0, t] ^p} \det(A_p(u) )^{-1/2}  \frac{1}{(2\pi)^{ p/2}}\det(A_p(u) )^{1/2}\\
 && \int_{\RR^p} \exp\left\{-\frac12 \xi^\top A_p(u) \xi\right\} \prod_{k=1}^p \left[ e^{-ix\xi_k}- 
e^{-iy\xi_k}\right]   d\xi du \,.  
\end{eqnarray*}  
Assume that $Z_1, \cdots, Z_p$ are jointly Gaussians with mean zero and covariance matrix 
$A_p ^{-1}(u) =(a_{ij}^{-1})_{1\le i,j\le p}$, where $A_p ^{-1}(u)  $ denotes the inverse matrix of $A_p(u)$ (recall the definition of $A_p$ in \ref{e.4.1}). Then we can write   for $x, y\in [-R, R]$
\begin{eqnarray}
\EE |\cL_t^x-\cL_t^y|^p   
&=&  \frac{1}{(2\pi)^{p/2}} \int_{[0, t] ^p} \det(A_p(u))^{-1/2}  \prod_{k=1}^p\EE 
      \left[ e^{-ix Z_k}- e^{-iy Z_k}\right]du\nonumber\ \\
 &\le &  \frac{1}{(2\pi)^{p/2}} \int_{[0, t] ^p} \det(A_p(u))^{-1/2}  \prod_{k=1}^p\left(\EE 
      \left| e^{-ix Z_k}- 
e^{-iy Z_k}\right|^{p}\right)^{1/p} du\nonumber\\
&\le &   C_{p, R} |x-y|^{\al p }  \int_{[0, t] ^p} \det(A_p(u))^{-1/2}  \prod_{k=1}^p\left(\EE 
      \left|   Z_k \right|^{\al p}\right)^{1/p} du   \,,  \label{e.37}
\end{eqnarray}  
where we used the fact that  for any $\al \in [0, 1]$, $|e^{i x}-e^{iy}|\le |x-y|= |x-y|^\al 
|x-y|^{1-\al}\le C_{\al, R} |x-y|^{\al}$
for any $x, y\in [-R, R]$.   By    identity  of expressing any moment via   the  variance
(see e.g. \cite[Equation (3.1.8)]{yaozhong2017analysis} or \eqref{e:p}) we have 
\begin{eqnarray}
 \left(\EE|Z_i|^{\al p}\right)^{1/p} \le C_{\al, p} (\EE|Z_i|^2 )^{\al   /2}\,. 
\end{eqnarray}
We want to bound $\EE|Z_i|^2$ appropriately.  Notice that $Z_i\sim N(0, a_{ii}^{-1})$,
where  $a_{ii}^{-1}$ is the $i$-th diagonal element of the inverse matrix of $A_p(u)$.  Hence,
\begin{eqnarray}
 \left(\EE|Z_i|^{\al p}\right)^{1/p} \le C_{\al, p} (a_{ii}^{-1}  )^{\al   /2}\,. 
\end{eqnarray}
Using the Cramer rule, we have 
\begin{eqnarray}
a_{ii}^{-1}=\frac{\det(\tilde A_{ii}(u))}{\det(A_p(u))}\,\,, 
\end{eqnarray}
where $\tilde A_{ii}(u)$ is the $(p-1)\times (p-1)$ matrix obtained from $A_p(u)$ by deleting the $i$-th row and 
the $i$-th column.  Thus by the second inequality of \eqref{e.4.2}, we have
\[
\det(\tilde A_{ii}(u))\le  u_1(u_2-u_1)\cdots (u_{i-1}-u_{i-2})(u_{i+1}-u_{i-1})\cdots (u_p-u_{p-1})\,.
\]
Combining this inequality with  the first inequality in \eqref{e.4.2}   for 
$\det(A_p(u))$ yields 
\begin{eqnarray*}
  a_{ii}^{-1} 
&  \le &
\left\{u_1(u_2-u_1)\cdots (u_p-u_{p-1})\exp\left[-2\al_*(u_p) u_p\right]\right\}^{-1} \\
&&\qquad\quad 
u_1(u_2-u_1)\cdots (u_{i-1}-u_{i-2})(u_{i+1}-u_{i-1})\cdots (u_p-u_{p-1})\\
&=& \exp\left[2\al_*(u_p) u_p\right] \frac{u_{i+1}-u_{i-1}}{(u_{i+1}-u_i)(u_{i}-u_{i-1})}\\
&=&  \exp\left[2\al_*(u_p) u_p\right] \left[\frac{1}{ u_{i+1}-u_i }+\frac{1}{ u_{i}-u_{i-1} }\right]\,.
\end{eqnarray*}

Thus, 
\[
\left(\EE|Z_i|^{\al p}\right)^{1/p} \le C_{p,\al, t}\left[ (u_{i+1}-u_i)^{-\al/2}+(u_{i  }-u_{i-1})^{-\al/2} \right]\,.
\]

Substituting the above inequality and the first inequality in \eqref{e.4.2} to 
\eqref{e.37}, we have, denoting $u_0=0$, 
\begin{eqnarray*}
\EE |\cL_t^x-\cL_t^y|^p  
&\le &   C_{\al, t, p, R} |x-y|^{\al p }  \int_{[0, t] ^p} \prod_{i=1}^p (u_i-u_{i-1})^{-1/2}\\
&& \left[ (u_{i+1}-u_i)^{-\al/2}+(u_{i  }-u_{i-1})^{-\al/2}\right]du \,. 
\end{eqnarray*}  
It  is easy to see that the above multiple integral is   bounded by a finite 
constant $C_{\al, p, t}$ for any $\al<1/2$.  This means  for any $\al<1/2$, and for any positive integer 
$p$,  we have 
\[  
\EE |\cL_t^x-\cL_t^y|^p  
 \le     C_{\al, t, p, R} |x-y|^{\al p} \,.  
\]
This proves the theorem by the Kolmogorov lemma 
(see \citet[Corollary 2.1]{yaozhong2017analysis}).  

 \end{pf} 
\begin{remark}\label{r.3.3} 
Ray (see e.g. \citet[Equation (2.214)]{marcusrosen}) 
used the first Ray-Knight theorem to give the following iterated logarithmic law 
for the local time of the Brownian motion:
\[
\limsup_{\de\rightarrow 0} \frac{|\cL_t^{x+\de}-\cL_t^x|}{\sqrt{\de \log \log 
\de}}=2\sqrt{\cL_t^x}\,,\quad \hbox{almost surely} 
\]
when $X_t$ is  the standard Brownian motion.  It seems hard to adopt the   bounds  in \eqref{e.37}
to show 
$|\cL_t^{y}-\cL_t^x|\le C_{R, t} \sqrt{|x-y||\log|x-y||}$.   One may need a more subtle bounds. 
\end{remark}

\section{Example}
\label{sec:5}
In this section, we conduct some numerical experiments to illustrate the convergence of $X$ by Monte Carlo simulations.

\begin{figure} 
\centering 
\subfigure[$\beta=0.8$]{\label{fig:subfig:a}
\includegraphics[width=0.8\textwidth]{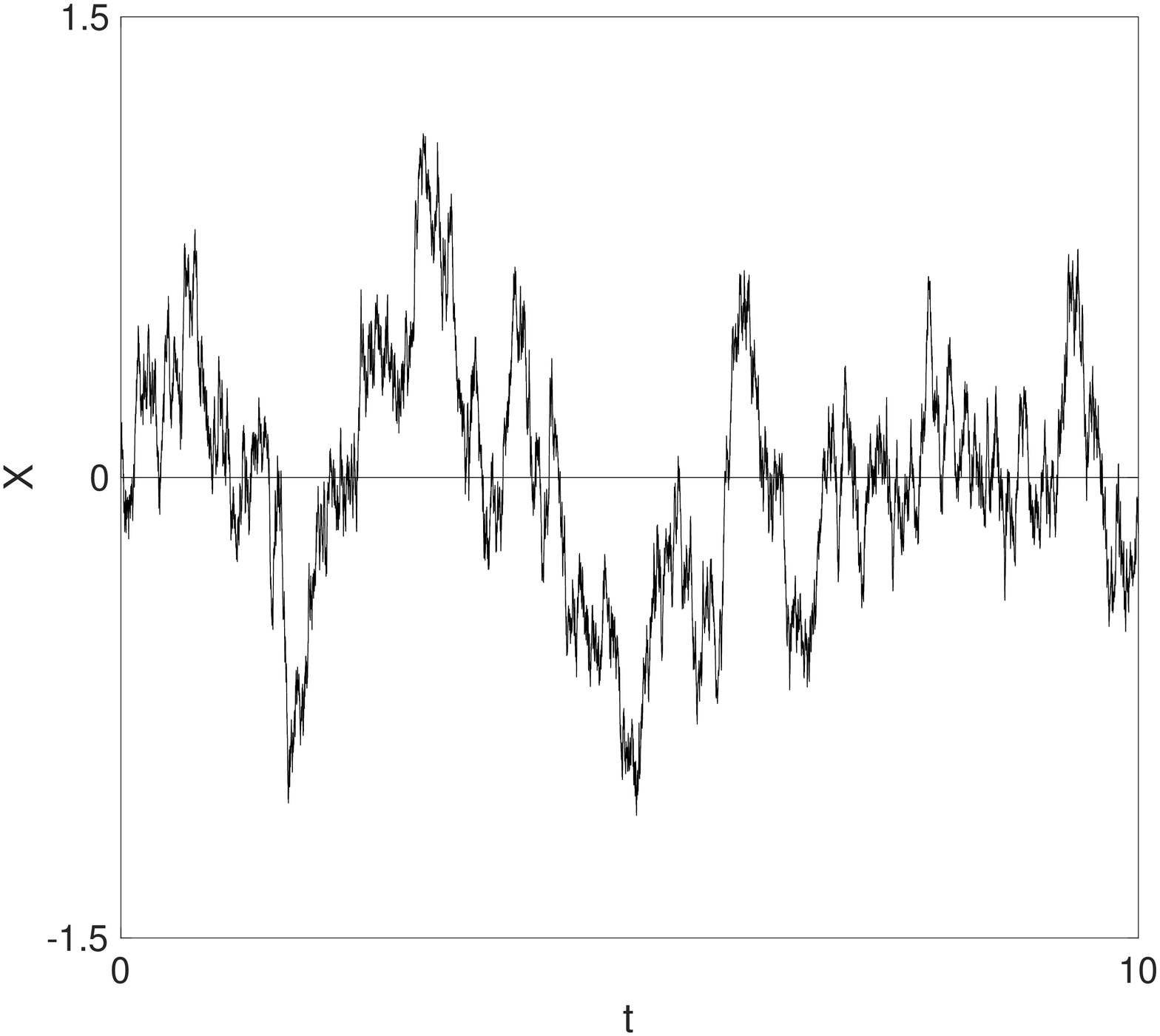}}
\subfigure[$\beta=2.0$]{\label{fig:subfig:b}
\includegraphics[width=0.8\textwidth]{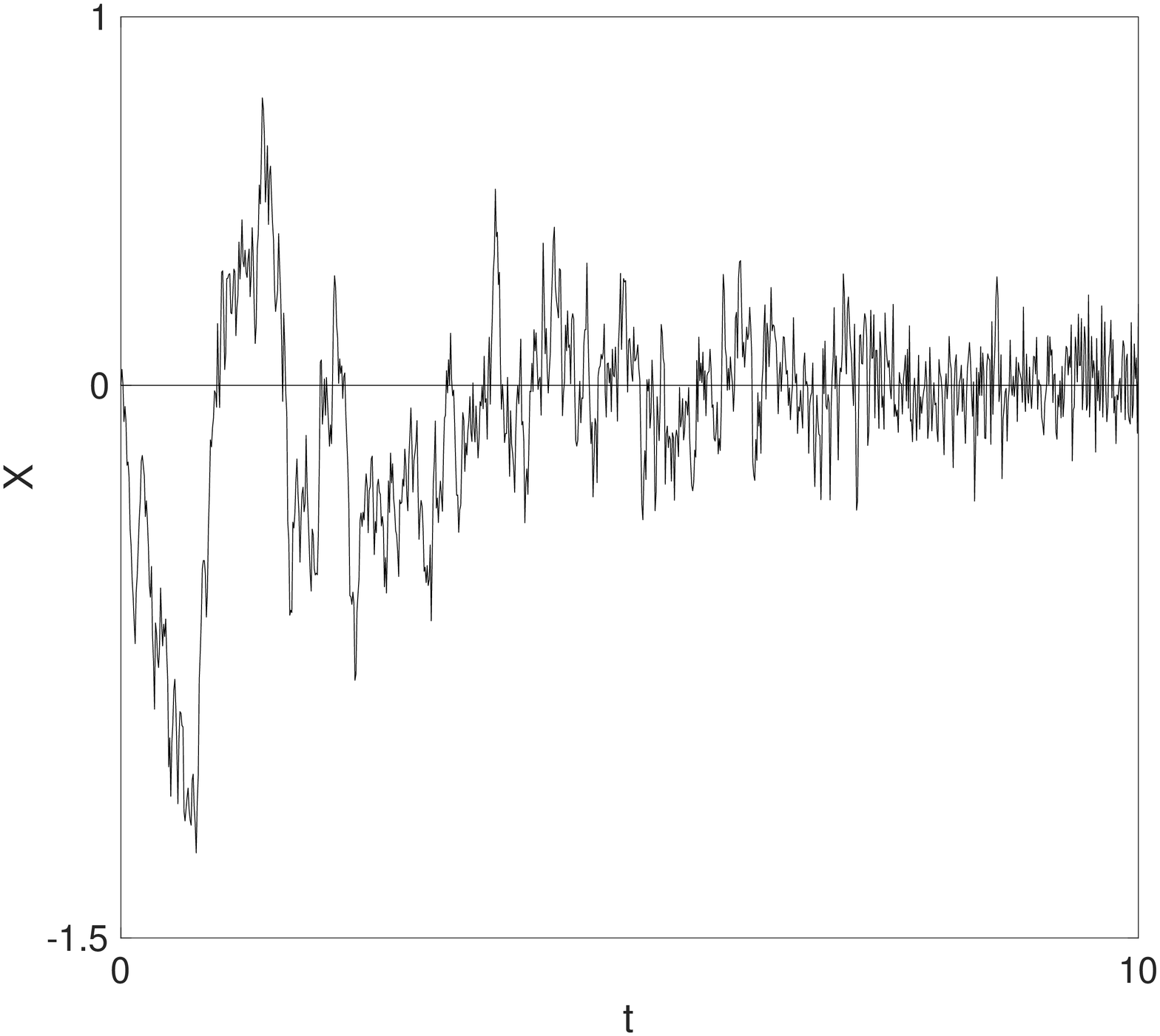}}
\caption{Simulation for infinite horizon Brownian bridges: $dX_t=-t^\beta X_t dt+ dW_t$. }
\label{fig:1Ex}
\end{figure}

We take  $\beta=0.8,2.0$ in Figure \ref{fig:1Ex} by simulating the following SDE:
\begin{equation}\label{e:1Ex}
dX_t=-t^\beta X_t dt +d W_t \;.
\end{equation}
 The parameters are set as time step $h=0.01$, initial value $X_0=0$. It can be seen that when the index $\beta$ is larger, the rate of convergence to zero is faster.

\begin{figure} 
\centering 
\subfigure[$\beta=0.5$]{\label{fig:subfig:a}
\includegraphics[width=0.8\linewidth]{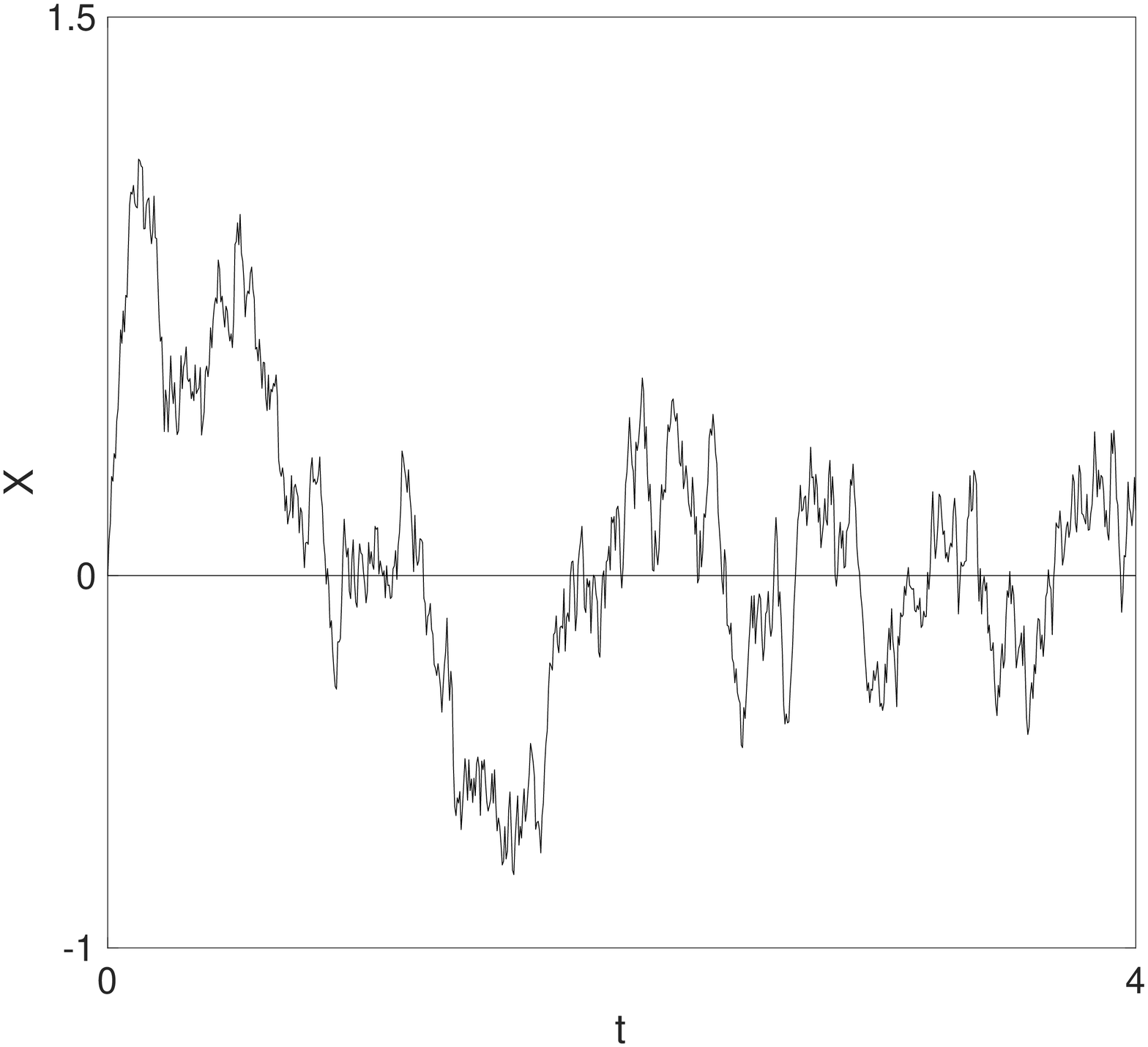}}
\subfigure[$\beta=1.5$]{\label{fig:subfig:b}
\includegraphics[width=0.8\linewidth]{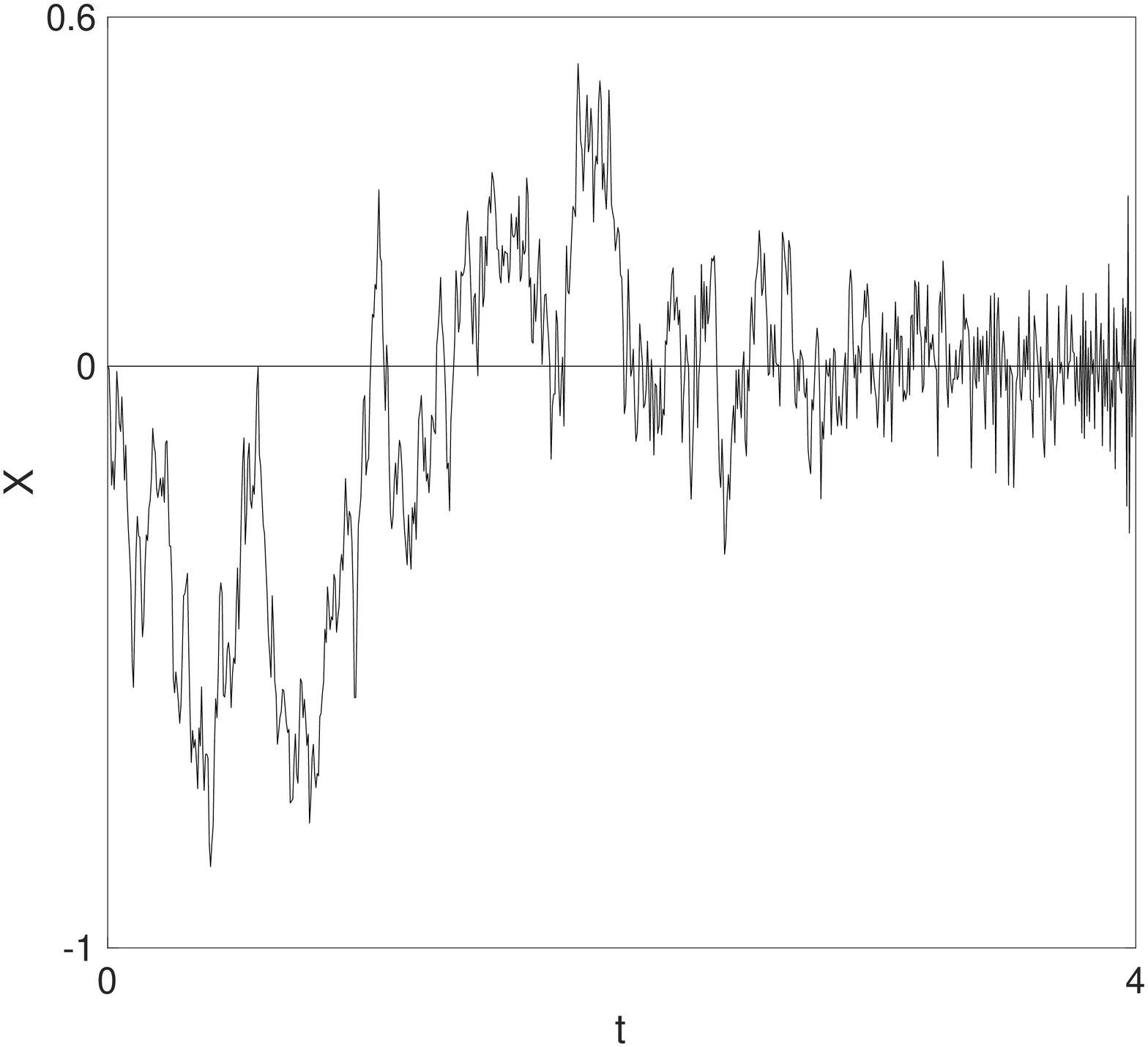}}
\caption{Simulation for infinite horizon Brownian bridges: $dX_t=-e^{\beta t} X_t dt+ dW_t$. }
\label{fig:2Ex}
\end{figure}

We take  $\beta=0.5,1.5$ in Figure \ref{fig:2Ex} by simulating the SDE:
\begin{equation}\label{e:2Ex}
dX_t=-e^{\beta t } X_t dt+ d W_t,
\end{equation}
The parameters are set as time step $h=0.005$, initial value $X_0=0$. The figure illustrates  that when the index $\beta$ is larger, the rate of convergence to zero is faster.

\section*{Acknowledgement}
 Y. Hu is supported by an NSERC discovery grant and a startup fund of University of Alberta.    Y. Xi is supported by the National Natural Science Foundation of China  (Grant No.  11631004, 71532001) and the China Scholarship Council.

\appendix
\section{Proofs of Lemmas \ref{e:Lop}-\ref{lem:2}}\label{A:A}
 Proof  of Lemma \ref{e:Lop}. 
It is easy to see that both the denominator and the numerator go to infinity when 
$x\rightarrow \infty$. We can use  the L'Hopital rule. The limit of the left hand side 
of \eqref{e.2.9} is the same as
\[
 \frac{g (x)  \exp\left(\kappa \int_0^x g(u)du\right)}{\kappa g(x)\exp\left( \kappa \int_0^x g(u)du\right)  }
 +\frac1{\kappa}\frac{g'(x)}{g(x)}\frac{\int_0^x\exp\left(\kappa \int_0^s g(u)du\right)ds}{ \exp\left( \kappa \int_0^x g(u)du\right)  }\,.
\]
The first summand is $1/\kappa$. Applying the  L'Hopital rule to the last fraction  of the second summand 
we see 
\[
\lim_{x\rightarrow \infty} \frac{\int_0^x\exp\left(\kappa \int_0^s g(u)du\right)ds}{ \exp\left( \kappa \int_0^x g(u)du\right)  }= \lim_{x\rightarrow \infty} \frac{ \exp\left(\kappa \int_0^x g(u)du\right) }{ \kappa g(x)\exp\left( \kappa \int_0^x g(u)du\right)  }=0\,. 
\]
This implies $\lim_{x\rightarrow \infty} \frac{g(x) \int_0^x\exp\left(\kappa \int_0^s g(u)du\right)ds}{\exp\left( \kappa \int_0^x g(u)du\right)  }=1/\kappa$. The lemma is then proved. 

\bigskip 
Proof of Lemma \ref{lem:1}. 
Without loss of generality, we assume that $0<t_1<t_2<\infty$.  Then
\begin{align}\nonumber
\sigma^2_{t_1, t_2}&=\mathbb E\left( \int_0^{t_2} \exp\left(-\int_s^{t_2} \alpha(u)du\right)dW_s -\int_0^{t_1}\exp\left(-\int_s^{t_1} \alpha(u)du\right)dW_s\right)^2 \\ \nonumber
&\le 2 \mathbb E\left[ \int_0^{t_1} \left[
\exp\left(-\int_s^{t_2} \alpha(u)du\right)- \exp\left(-\int_s^{t_1} \alpha(u)du\right)
\right] dW_s\right]^2\\\nonumber
& + 2\mathbb E\left(\int_{t_1}^{t_2} \exp\left(-\int_s^{t_2} \alpha(u)du\right)dW_s\right)^2 \\ \label{ine}
&= 2   \int_0^{t_1} \left[
\exp\left(-\int_s^{t_2} \alpha(u)du\right)- \exp\left(-\int_s^{t_1} \alpha(u)du\right)
\right]^2 ds \nonumber\\
&\qquad\qquad  + 2 \int_{t_1}^{t_2} \exp\left(-2\int_s^{t_2} \alpha(u)du\right)ds=:I_1+I_2 .
  \end{align}
Denote 
\begin{equation*}
A_*(t)=\int_0^t \al(u) du\,, \quad t\ge 0\,. 
\end{equation*}
We estimate the first term in  \eqref{ine} as  follows. For any $0\le\gamma\le 1$, we have 
\begin{align}
I_1&= 2 \left(1-\exp\left(-\int_{t_1}^{t_2}\alpha(u)du\right)\right)^2\int_0^{t_1}\exp\left( -2\int_s^{t_1} \al(u) du\right)ds
\nonumber\\
&\le C_\ga \left(\int_{t_1}^{t_2}\alpha(u)du\right)^{2\ga}  \int_0^{t_1}\exp\left[  -2A_*(t_1)+2A_*(s) \right]ds \nonumber\\  
&\le C_\ga \left(\int_{t_1}^{t_2}\alpha^*(u)du\right)^{2\ga}   \frac{\int_0^{t_1}\exp\left[  2A_*(s) \right]ds}{ \exp\left[   2A_*(t_1) \right] }\nonumber\\ 
&\le C_{\ga} 
 \left(t_2-t_1\right)^{2\ga} \left(\al^*(t_2)\right)^{2\ga}    (\al(t_1))^{-1}\,, 
 \label{e:1term} 
\end{align}
where in the first inequality  we have used the inequality that for any $\ga\in [0, 1]$,
$1-e^{-x}\le C_\ga x^\ga\,, x\ge 0$ 
 for some constant $C_\ga$, and the last inequality follows from Lemma \ref{e:Lop}.
Now we estimate the second term in \eqref{ine}. Since $\al(t)$ is a positive function,
we have 
\begin{align*}
\int_{t_1}^{t_2}\exp\left( -2\int_s^{t_2}\alpha(u)du \right)ds
\le t_2-t_1\,. 
\end{align*}
On the other hand, we have 
\begin{align*}
\int_{t_1}^{t_2}\exp\left( -2\int_s^{t_2}\alpha(u)du \right)ds
&\le  \int_0^{t_2}\exp\left( -2\int_s^{t_2}\alpha(u)du \right)ds
 \le  C (\al(t_2))^{-1} \,, 
\end{align*}
where in the above last inequality we use Lemma \ref{e:Lop} again as in the above 
argument for $I_1$.  
Therefore,  
\begin{equation}\label{e:2term}
\int_{t_1}^{t_2}\exp\left( -2\int_s^{t_2}\alpha(u)du \right)ds
 \le    C (t_2-t_1)\wedge (\al(t_2))^{-1}\,. 
\end{equation}
Combining  \eqref{e:1term} and \eqref{e:2term}, for $0\le\ga\le1/2$, we have 
\begin{align*}
\sigma^2_{t_1, t_2}&:=\mathbb E|X_{t_2}-X_{t_1}|^2\leq C_{\ga} 
 \left|t_2-t_1\right|^{2\ga} \left(\al^*(t_1\vee t_2)\right)^{2\ga}    (\al(t_1\wedge t_2))^{-1}\\
 &+   C|t_2-t_1|\wedge (\al(t_1\vee t_2))^{-1}\\
 &\le C_{ \ga}|t_2-t_1|^{2\ga}\left[ \left(\al^*(t_1\vee t_2)\right)^{2\ga}   (\al(t_1\wedge t_2))^{-1} +
(\al(t_1\vee t_2))^{-(1-2\ga)} \right] \,,
\end{align*}
where we used the inequality 
\[
a\wedge b= (a\wedge b)^\ga \cdot 
(a\wedge b)^{1-\ga}\le a^\ga b^{1-\ga}, 0\le\ga\le 1\,.
\]
Thus,  we  have proved the lemma.

\bigskip 
Proof of Lemma \ref{lem:2}. 
Let  $m, n$ be integers, $m\geq2$, $n\geq1$. Since $X_{t_2}-X_{t_1}$ is a one-dimensional Gaussian process with mean zero and variance $\sigma^2_{t_1,t_2}$ 
(defined by \eqref{e.2.6})  we can express its    moments by this variance
(see e.g. \cite[Equation (3.1.8)]{yaozhong2017analysis})
\begin{equation}
\mathbb E|X_{t_2}-X_{t_1}|^m=\left\{
\begin{array}{lcl}\label{e:p}
 \frac{(2n)!(\sigma^2_{t_1,t_2})^{n}}{2^nn!},& & \text{if } m=2n \text{ is even},\\
0, & & {\text{if } m=2n+1 \text{ is odd}}\,. 
\end{array} \right.
\end{equation}
From now on we assume  $m$ is an even integer. Denote 
\begin{equation*}
\rho(x;k):= C_{\ga}x^{\ga} 
 k ^{-\beta(1-2\ga)/2}  \,.
\end{equation*}
From condition (i) in Theorem \ref{th:1}, we have $\al(t)\ge Ct^{\beta}$ for some $C>0$, which yields
\begin{align*}
\sigma_{t_1,t_2}^2&\le  C_{ \ga}|{t_2}-{t_1}|^{2\ga}\left[ \left(\al^*(k+1)\right)^{2\ga}   (\al(t_1\wedge t_2))^{-1} +
(\al(t_1\vee t_2))^{-(1-2\ga)} \right]\\
&\le  C_{ \ga}|{t_2}-{t_1}|^{2\ga}\left[ \left(\al^*(t_1\wedge t_2+1)\right)^{2\ga}   (\al(t_1\wedge t_2))^{-1} +
(\al(t_1\wedge t_2))^{-(1-2\ga)} \right]\\
&\le C_\ga|{t_2}-{t_1}|^{2\ga}\left[  C   (t_1\wedge t_2)^{-\beta} +
 C (t_1\wedge t_2) ^{-\beta(1-2\ga)} \right]\\
&\le C_\ga|{t_2}-{t_1}|^{2\ga}\left[  C   k^{-\beta} +
 C k ^{-\beta(1-2\ga)} \right]\\
&\le \left(\rho(|t_2-t_1|;k)\right)^2\,.
\end{align*}
 Lemma \ref{lem:1} and \eqref{e:p} imply that
\begin{equation}\label{e:inp}
\mathbb E|X_{t_2}-X_{t_1}|^m\le \frac{m!\sigma_{t_1t_2}^m}{2^{m/2}(m/2)!}\le \frac{m!(\rho(|t_2-t_1|;k))^{m}}{2^{m/2}(m/2)!}\;.
\end{equation}
Set 
\begin{equation}
B_k:=\int_k^{k+1}\int_k^{k+1}\frac{|X_t-X_s|^m}{(\rho(|t-s|;k))^{m}}dsdt\,. 
\end{equation}
Then $B_k$ is finite. Take $\Psi(x)=x^m$. The inequality \eqref{e:inp}  implies 
\begin{equation}
\mathbb E(B_k)=\mathbb E\int_k^{k+1}\int_k^{k+1}\Psi\left(\frac{|X_t-X_s|}{\rho(|t-s|;k)}\right)dsdt\le \frac{m!}{2^{m/2}(m/2)!}\;.
\end{equation}
For any $N>1$, we have
\begin{equation*}
\mathbb E\left( \sum_{k=1}^\infty\frac{B_k}{k^N}\right)=\sum_{k=1}^\infty\frac{\mathbb E(B_k)}{k^N}< \infty.
\end{equation*}
This implies that 
\begin{equation}
  R_{N,m}:=\sum_{k=1}^\infty\frac{B_k}{k^N} ~\textup{ is an   almost surely finite  random constant }  \,.
\end{equation}
Since all $B_k$  are positive, we see 
\begin{equation}
B_k\le R_{N,m}k^N ~\textup{ for all positive number }  N>1 \textup{ and all integer } k\ge 1.
\end{equation}
By virtue of the Garsia-Rodemich-Rumsey inequality, see e.g., \citet[Theorem 2.1]{yaozhong2017analysis}, we can choose $m\ga>2$, $2N/m<\beta$, and $2\ga<1-2N/(\beta m)$  such that 
for any $t_1, t_2\in [k, k+1]$
\begin{align*}
|X_{t_2}-X_{t_1}|&\le 8\int_0^{|t_2-t_1|}\Psi^{-1}\left( \frac{4B_k}{u^2}\right)\rho^{\prime}(u;k)du\\
&=8\left(4B_k\right)^{1/m}\ga C_\ga  
 k ^{-\beta(1-2\ga)/2} \int_0^{|t_2-t_1|} u^{\ga-1-\frac{2}{m}}du\\
&\le \frac{8\left(4R_{N,m}\right)^{1/m}\ga}{\ga-\frac{2}{m}}
C_\ga  k^{-\beta\left[1-2\ga-2N/(\beta m)\right]/2}  
\end{align*}
since $|t_2-t_1|\le 1$,  where $\Psi^{-1}(\cdot)$ is the inverse function of $\Psi(\cdot)$. This proves
 the lemma. 

\section{Proof of Lemma \ref{lem:5}}\label{B}
 First let us recall a well-known result for Gaussian random variables. 
Let  $Z_1, \cdots, Z_m$   be  
a  set of  centered jointly Gaussian  random  variables with covariance matrix 
$F=\left(\EE(Z_iZ_j)\right)_{1\le  i,j\le p}$.   It is elementary from the  well known form of  the  multivariate normal  distribution that (cf. \citet{berman1973local})  the determinant of $F$ has the following representation: 
\begin{equation}
\det(F)=\var(Z_1)\var(Z_2|Z_1)\cdots\var(Z_p|Z_1, \cdots, Z_{p-1})\,,
\label{e.4.3} 
\end{equation}
where    
\begin{eqnarray*}
\var(Z|Y_1, \cdots,Y_k)
&=&\EE\left\{\left[Z-\EE(Z|Y_1, \cdots,Y_k)\right]^2 |Y_1, \cdots,Y_k)\right\}\\
&=&\EE\left\{\left[Z-\EE(Z|Y_1, \cdots,Y_k)\right]^2  \right\} 
\end{eqnarray*} 
denotes  the conditional variance of $Z$ given $Y_1, \cdots, Y_k$  and the above last identity follows from the fact that  $Z-\EE(Z|Y_1, \cdots,Y_k)$ is independent of $  Y_1, \cdots,Y_k $.   For the solution of \eqref{e:FBB}, we have
\begin{eqnarray}
\var(X_t|X_s)
&=&\EE\left[(X_t-\EE(X_t|X_s))^2\right]
=\EE\left[\left(\int_s^t e^{-\int_r^t \al(u) du}dW_r\right)^2\right]\nonumber\\
&= & \int_s^t e^{-2\int_r^t \al(u) du }d r\ge (t-s)\exp\left(-2\al^*(t)(t-s)\right)\,. \label{e.4.4} 
\end{eqnarray}
For $u_1<u_2<\cdots<u_p$ applying the   identity \eqref{e.4.3} to $X_{u_1}\,, \cdots, X_{u_p}$ and using the above  inequality we have 
\begin{eqnarray}
\det(A)&=&\mathrm{Var}(X_{u_1}) \var(X_{u_2}|X_{u_1})\cdots \var(X_{u_p}|
X_{u_1}, \cdots, X_{u_{p-1}})\nonumber\\ 
&\ge & u_1(u_2-u_1)\cdots (u_p-u_{p-1})\exp\left\{-2\al^*(u_p) u_p\right\}\,. \label{e.3.24} 
\end{eqnarray}
This proves the first   inequality in \eqref{e.4.2}. 
On the other hand,   by \eqref{e.4.4} we  have for any $0\le s <t<\infty$, 
\begin{equation*}
\var(X_t|\cF_s) \le t-s \,. 
\end{equation*}  
This can be used to prove the second  inequality in \eqref{e.4.2}.

\bibliography{mybibfile}

\end{document}